\documentclass{amsproc}
\usepackage{graphicx}
\usepackage{amssymb}
\usepackage{epstopdf}
\usepackage{amsmath,amsthm,amscd,amssymb}
\usepackage{latexsym}
\usepackage[colorlinks,citecolor=red,pagebackref,hypertexnames=false]{hyperref}
\usepackage{geometry}                
\numberwithin{equation}{section}

\theoremstyle{plain}
\newtheorem{theorem}{Theorem}[section]
\newtheorem{lemma}[theorem]{Lemma}

\newtheorem{conjecture}[theorem]{Conjecture}

\theoremstyle{definition}
\newtheorem{definition}[theorem]{Definition}

\newtheorem{case[theorem]}{Case}

\theoremstyle{remark}
\newtheorem{remark}[theorem]{Remark}

\numberwithin{equation}{section}

\begin{document}

\title{\parbox{14cm}{\centering{On sets of directions determined by subsets of ${\Bbb R}^d$}}}


\author{Alex Iosevich, Mihalis Mourgoglou, and Steven Senger}

\begin{abstract} Given $E \subset \mathbb{R}^d$, $d \ge 2$, define ${\mathcal D}(E) \equiv \left\{ \frac{x-y}{|x-y|}: x,y \in E \right\} \subset S^{d-1},$ to be the set of directions determined by $E$. We prove that if the Hausdorff dimension of $E$ is greater than $d-1$, then $\sigma({\mathcal D}(E))>0$, where $\sigma$ denotes the surface measure on $S^{d-1}$. In the process, we prove some tight upper and lower bound for the maximal function associated with the Radon-Nikodym derivative of the natural measure on ${\mathcal D}$. This result is sharp since the conclusion fails to hold if $E$ is a $(d-1)$-dimensional hyper-plane. This result can be viewed as a continuous analog of a recent of Pach, Pinchasi, and Sharir (\cite{PPS04}, \cite{PPS07}) on directions determined by finite subsets of $\mathbb{R}^d$. We also discuss the case when the Hausdorff dimension of $E$ is precisely $d-1$, where some interesting counter-examples were previously obtained by Simon and Solomyak (\cite{SiSo06}) in the planar case. In the case when the Hausdorff dimension of $E$ equals $d-1$ and set is rectifiable and not contained in a hyper-pane, T. Orponen and T. Sahlsten (\cite{OS10}) recently proved that the Lebesgue measure of the set of directions is still positive, in response to the conjecture stated in this paper. 

At the end of this paper we show that our continuous results can be used to recover and in some case improve the exponents for the corresponding results in the discrete setting for large classes of finite point sets. In particular, we prove that a finite point set $P \subset {\Bbb R}^d$, $d \ge 3$, satisfying a certain discrete energy condition (Definition \ref{adaptablemama}), determines $\gtrapprox \# P$ distinct directions. 

\end{abstract} 

\maketitle


\section{Introduction}

\vskip.125in 

A large class of Erd\H os type problem in geometric combinatorics asks whether a large set of points in Euclidean space determines a suitably large sets of geometric relations or objects. For example, the classical Erd\H os distance problem asks whether $N$ points in ${\Bbb R}^d$, $d \ge 2$, determines $\gtrapprox N^{\frac{2}{d}}$ distinct distances, where here, and throughout, $X \lessapprox Y$, with the controlling parameter $N$ means that for every $\epsilon>0$ there exists $C_{\epsilon}>0$ such that $X \leq C_{\epsilon} N^{\epsilon} Y$. See, for example \cite{BMP00}, \cite{Ma02}, \cite{P05}, \cite{PS04}, \cite{Sz97} and the references contained therein for a thorough description of these types of problems and recent results. 

Continuous variants of Erd\H os type geometric problems have also received much attention in recent decades. Perhaps the best known of these is the Falconer distance problem, which asks whether the Lebesgue measure of the distance set $\{|x-y|: x,y \in E\}$ is positive, provided that the Hausdorff dimension of $E \subset {\Bbb R}^d$, $d \ge 2$, is greater than $\frac{d}{2}$. See \cite{Erd05} and \cite{W99} for the best currently known results on this problem. See also \cite{B94}, \cite{IS10}, \cite{Mat87} and \cite{M95}. Also see \cite{CEHIT10} for the closely related problem on finite point configurations. For related problems under the assumption of positive Lebesgue density, see, for example, \cite{B86}, \cite{FKW90}, and \cite{Z06}.

In this paper we study the sets of directions determined by subsets of the Euclidean space. In the discrete setting, the problem of directions was studied in recent years by Pach, Pinchasi, and Sharir. See \cite{PPS04} and \cite{PPS07}. In the latter paper they prove that if $P$ is a set of $n$ points in ${\Bbb R}^3$, not all in a common line or plane, then the pairs of points of $P$ determine at least $2n-5$ distinct directions if $n$ is odd, and at least $2n-7$ distinct directions if $n$ is even. Our main result can be viewed as a continuous variant of this result where finite point sets are replaced by infinite sets of a given Hausdorff dimension. An explicit quantitative connection between our main result on directions (Theorem \ref{main} below) and the work of Pach, Pinchasi, and Sharir is made in Section \ref{discrete} below. We show that a finite set $P$, satisfying the $(d-1+\epsilon)$-adaptability assumption (see Definition \ref{adaptablemama} below), determines $\gtrapprox \# P$ distinct directions. In dimensions two and three, this result is weaker than then the result of Pach, Pinchasi and Sharir described above. However, in dimensions four and higher, our result gives, to the best of our knowledge, the only known bounds. 

In the finite field setting, the problem of directions was previously studied by the first listed author, Hannah Morgan and Jonathan Pakianathan. See \cite{IM10} and the references contained therein. 

\begin{definition} Given $E \subset {\Bbb R}^d$, $d \ge 2$, define 

$$ {\mathcal D}(E)=\left\{ \frac{x-y}{|x-y|}: x,y \in E, \text{ and } x \not=y \right\} \subset S^{d-1},$$ the set of directions determined by $E$. 

\end{definition} 

Our main results are the following. 

\begin{theorem} \label{main} Let $E \subset {\Bbb R}^d$, $d \ge 2$, of Hausdorff dimension greater than $d-1$. Let $\nu$ denote the probability measure on ${\mathcal D}(E)$ given by the relation 
$$ \int h(\omega) d\nu(\omega)=\int \int h \left( \frac{x-y}{|x-y|} \right) d\mu(x) d\mu(y),$$ where $\mu$ is a Frostman measure on $E$, and let $\nu^{\Theta}$ denote the same measure corresponding to the set $E^{\Theta}=\{ \Theta x: x \in E\}$, where $\Theta \in O(d)$, the orthogonal group. Let $\nu_{\epsilon}^{\Theta}(\omega)=\epsilon^{-(d-1)} \nu^{\Theta}(B(\omega, \epsilon))$, where $B(\omega, \epsilon)$ is the ball of radius $\epsilon$ centered at $\omega \in S^{d-1}$. Then 
$$ \nu_{\epsilon}^{\Theta}(\omega)=M^{\Theta}(\omega)+R^{\Theta}_{\epsilon}(\omega)$$ where 
$$ \int \sup_{0<\epsilon<c} |R^{\Theta}_{\epsilon}(\omega)| d\omega \lesssim c^{s-(d-1)},$$ 
$$ \int M^{\Theta}(\omega) d\omega \lesssim 1,$$ and there exists $\Theta_0 \in O(d)$ such that 
$$ \int M^{\Theta_0}(\omega) d\omega \gtrsim 1.$$ 

In particular, 
\begin{equation} \label{orgasm} \sigma({\mathcal D}(E))>0, \end{equation} 

\vskip.125in

where $\sigma$ denotes the Lebesgue measure on $S^{d-1}$.
\end{theorem}

\vskip.125in 

\begin{remark} Pertti Mattila recently pointed out to us that (\ref{orgasm}) follows from Theorem 10.11 in \cite{M95}. Our method, which uses Fourier analysis, allows us to obtain more detailed information about the direction set measure and the associated Randon-Nikodym derivative. \end{remark} 

\begin{remark} The approaches to similar problems in geometric measure theory (see e.g. \cite{Fal86}, \cite{Erd05}, \cite{W99}) typically involve constructing a measure on a set under consideration, (in this case-- directions), and then proving, using Fourier transform methods, that this measure is in $L^2$ or in $L^{\infty}$. While our approach is also Fourier based, we prove that the measure of an $\epsilon$-ball centered at any point in a subset of the the direction set, $\mathcal{F}$, of positive Lebesgue measure equals a quantity comparable to the Lebesgue measure of that ball plus an error $O(\epsilon^{s})$. Provided that the Hausdorff dimension of the underlying set $E$ is greater than $d-1$ this shows that the restriction of our measure to $\mathcal{F}$ is absolutely continuous with respect to the Lebesgue measure with Radon-Nikodym derivative bounded from above and below by positive constants, which implies that every Borel subset of $\mathcal{F}$ is contained in ${\mathcal D}(E)$ and allows us to conclude that the set of directions has positive Lebesgue measure on the sphere. This approach is quite reminiscent of the techniques used to study geometric combinatorics problems in the finite field setting. See, for example, (\cite{HIKR10}) and the references contained therein. 
\end{remark} 

\begin{remark} It is not difficult to check that if $E$ is a $(d-1)$-dimensional Lipschitz surface in ${\Bbb R}^d$, which is not contained in a $(d-1)$-dimensional plane, then $\sigma({\mathcal D}(E))>0$. It is reasonable to conjecture that the same conclusion holds if $E$ is merely a $(d-1)$-dimensional rectifiable subset of ${\Bbb R}^d$. We discuss the purely non-recitifiable case in the Subsection \ref{sharpnessofresults} below. \end{remark} 

\begin{remark} It is interesting to contrast this result with the Besicovitch-Kakeya conjecture (see e.g. \cite{W99II} and the references contained therein), which says that any subset of ${\Bbb R}^d$, containing a unit line segment in every direction has Hausdorff dimension $d$. On the other hand, Theorem \ref{main} says that Hausdorff dimension greater than $d-1$ is sufficient for the set to contain endpoints of a segment of some length pointing in the direction of a positive proportion of vectors in $S^{d-1}$. 
\end{remark} 

\subsection{Sharpness of the main results:} \label{sharpnessofresults} 

Theorem \ref{main} cannot be improved in the following sense. Suppose that $E$ is contained in a $(d-1)$-dimensional hyper-plane. Then $\sigma({\mathcal D}(E))=0$. It follows that the conclusion of Theorem \ref{main} does not in general hold if the Hausdorff dimension of $E$ is less than or equal to $d-1$. 

Another very different sharpness example comes from the theory of distance sets. Let $E_q$ denote the $q^{-\frac{d}{s}}$-neighborhood of 

$$q^{-1} \left( {\Bbb Z}^d \cap {[0,q]}^d \right),$$ where ${\Bbb Z}^d$ denote the standard integer lattice and $0<s<d$. It is known that if $q_i$ is a sequence of integers given by $q_1=2$, $q_{i+1}>q_i^i$, then the Hausdorff dimension of 

$$E= \cap_i E_{q_i}$$ is $s$. See, for example, \cite{Fal86}, \cite{Falc86}. Observe that 

$$ \sigma({\mathcal D}(E_q)) \approx q^{-\frac{d(d-1)}{s}} \cdot q^d$$ since the number of lattice points in ${[0,q]}^d$, $d \ge 2$, with relatively prime coordinates is equal to 

$$\frac{q^d}{\zeta(d)}(1+o(1)),$$ where $\zeta(t)$ is the Riemann zeta function. See, for example, \cite{Landau27}. It follows that 

$$ \sigma({\mathcal D}(E_q)) \to 0 \ \text{as} \ q \to \infty$$ if $s<d-1$. It follows that $\sigma({\mathcal D}(E))=0$. This example does not rule out $s=d-1$ and one might reasonably conjecture, consistent in spirit with the result due to Pach, Pinchasi, and Sharir stated above, that if the Hausdorff dimension of $E$ is equal to $d-1$, then (\ref{orgasm}) holds if and only $E$ is not a subset of a single $(d-1)$-dimensional hyper-plane. This, however, is not true. A result due to Simon and Solomyak (\cite{SiSo06} shows that for every self-similar set of Hausdorff dimension one satisfying an additional mild condition, the Lebesgue measure of ${\mathcal D}(E)$ is zero. In particular, if $E$ is the four-cornered Cantor set known as the Garnett set (see e.g. \cite{GV03}), then the Hausdorff dimension of $E$ is $one$ and the Lebesgue measure of ${\mathcal D}(E)$ is zero. It is not difficult to use Simon and Solomyak's result to construct a set $E$ of Hausdorff dimension $d-1$ in ${\Bbb R}^d$ that is not contained in a hyperplane and the Lebesgue measure of ${\mathcal D}(E)$ is zero. 

\vskip.125in 

In the realm of rectifiable sets, we believe that Theorem \ref{main} can be strengthened as follows. 

\begin{conjecture} \label{whatthefuck} Let $E \subset {\Bbb R}^d$, $d \ge 2$, of Hausdorff dimension $d-1$. Suppose that $E$ is rectifiable and is not contained in a hyper-plane. Then (\ref{orgasm}) holds. 
\end{conjecture} 

\begin{remark} After this paper was submitted, this conjectured was resolved via a very nice argument due to Orponen and Sahlsten (\cite{OS10}). \end{remark}

\subsection{Structure of the paper} Theorem \ref{main} is proved in Section \ref{mainsection} below.  In Section \ref{discrete} we describe an explicit connection between the main results of the paper and the discrete problems, such as those studied by Pach, Pinchasi, and Sharir. 

\subsection{Acknowledgments} The authors are deeply grateful to Pertti Mattila for several very helpful remarks that significantly improved this paper. 

\vskip.125in 

\section{Proof of Theorem \ref{main}}
\label{mainsection}
\vskip.125in
Let $s$ be the Hausdorff dimension of $E$. Now, although the set $\mathcal{D}(E)$ is a subset of the $(d-1)$-dimensional sphere, in the arguments to follow, it is convenient to work with sets of slopes of line segments defined by pairs of points in the set $E$. If $p$ and $q$ are two points, with coordinates $(p_1, p_2, ..., p_d)$ and $(q_1, q_2, ..., q_d)$, then we define the slope of the line segment between $p$ and $q$ as the $(d-1)$-tuple
$$\left\lbrace \frac{p_1-q_1}{p_d-q_d}, \frac{p_2-q_2}{p_d-q_d}, ...,
\frac{p_{d-1}-q_{d-1}}{p_d-q_d} \right\rbrace.$$
\vskip.125in
It is not difficult to see that if the $(d-1)$-dimensional Lebesgue measure of the set of slopes determined by $E$ is positive, then (\ref{orgasm}) holds. With a slight abuse of notation, we will refer to the set of slopes as $\mathcal{D}(E)$ as well.
It is convenient to extract two subsets from $E$, separated from each other in at least one of the coordinates. To this end, we have the following construction.
\begin{lemma} \label{twosets} If $\mu$ is a Frostman probability measure on $E \subset {\Bbb R}^d,$ with Hausdorff dimension greater than $d-1$, then there exist $c_1, c_2$ positive constants and $E_1, E_2$ subsets of $E$ such that $\mu(E_j) \geq c_1 >0,$ for $j=1,2$ and
\begin{equation*}
\max\limits_{1 \leq k \leq d} (\inf \{|x_k-y_k|: x\in E_1, y \in E_2\})\geq c_2 > 0.
\end{equation*}
\end{lemma}
\vskip.125in

We will employ a stopping time argument. Define $C_0$ to be the constant in the Frostman condition, 

$$\mu(B_r)\leq C_0 r^{s}.$$ 

\vskip.125in

Let $[0,1]^d$ be the unit cube in $\Bbb{R}^d$, and subdivide it into $4^d$ smaller cubes of side-length $\frac{1}{4}.$ Choose $2^d$ collections of $2^d$ sub-cubes each, such that no two cubes of the same collection touch each other. Then by the pigeon-hole principle, at least one of them has measure greater than or equal to $\frac{1}{2^d}.$ If there are two cubes, $Q_1$, and $Q_1'$, in the same collection, such that $\mu(Q_1), \mu(Q_1') \geq \frac{c}{2^d}$ for some $c>0,$ then we are done. If not, there exists a cube $Q_1$ of side-length $\frac{1}{4}$, so that $\mu(Q_1) \geq \frac{1}{2^d}.$ Then we repeat the same procedure on the cube $Q_1$. Now, either we have two cubes, $Q_2$ and $Q_2'$, with $\mu(Q_2)$, and $\mu(Q_2') \geq \frac{c}{2^{2d}}$, for some $c>0$, which are in the same collection, or we do not. If we do not, then again, there must be a cube, $Q_2$, with side length $\frac{1}{4^2}$, so that $\mu(Q_2) \geq \frac{1}{2^{2d}}$. We can repeat this process, and at each stage check for two cubes, from the same collection, with the requisite measure. Let the integer $n$ depend on $C_0$. If we fail to find two such cubes at the n-th iteration, we obtain a cube $Q_n$ of side-length $\frac{1}{4^n}$ for which $\mu(Q) \geq \frac{1}{2^{dn}}.$ By the Frostman measure condition, there exists $C_0 > 1$ such that
\begin{equation*}
\frac{1}{2^{dn}} \leq \mu(Q) \leq C_0 \frac{1}{4^{sn}},
\end{equation*} which is true if $n \leq \frac{\log_2(C_0)}{(2s - d)}.$ So, picking $n > \frac{\log_2(C_0)}{(2s - d)},$ it is only true whenever $s < \frac{d}{2},$ and since $s > d-1 $, we have a contradiction.

\begin{figure}
\centering
\includegraphics[scale=1]{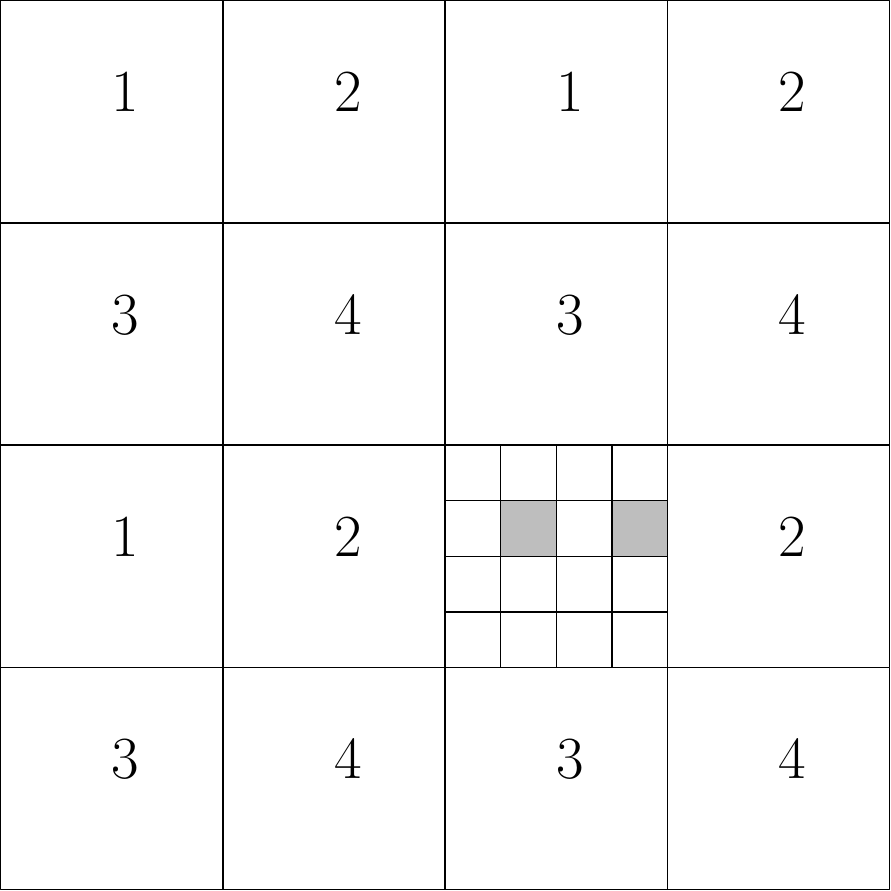}
\caption{The first decomposition into four collections of four cubes each is shown with a 1 in the cubes of the first collection, a 2 in the cubes of the second collection, etc... In this case, the first decomposition was not enough, and a positive proportion of the mass was in the lower-right cube of the first collection. After the second iteration, there are two shaded boxes, representing $E_1$ and $E_2$.}
\label{cubeFigure}
\end{figure}

\vskip.125in
Let $x$ and $y$ be points in $\mathbb{R}^d$ with coordinates $(x_1,
x_2, ..., x_d)$ and $(y_1, y_2, ..., y_d)$. Apply Lemma \ref{twosets} to $E$. Without loss of generality, let the sets $E_1$ and $E_2$ be separated in the $d$-th coordinate. Let $\mu_1$ and $\mu_2$ be restrictions of $\mu$ to the sets $E_1$ and $E_2$, respectively. Let $t=(t_1, t_2, \dots, t_{d-1})$. For slopes $t \in [\frac{1}{2},1]^{d-1}$, define $\nu_\epsilon(t)$ to be the quantity 

$$\frac{1}{\epsilon^{d-1}}\mu \times \mu \left\lbrace (x,y) \in E_1\times E_2 : t_1-\epsilon \leq \frac{x_1-y_1}{x_d-y_d} \leq t_1 + \epsilon, ..., t_{d-1}-\epsilon \leq \frac{x_{d-1}-y_{d-1}}{x_d-y_d} \leq t_{d-1} + \epsilon \right\rbrace.$$

Since $x_d-y_d$ is guaranteed to be more than $c_2$ by Lemma \ref{twosets}, we can multiply each inequality through by the denominator to get that
\begin{align*}
\nu_\epsilon(t) &\approx \frac{1}{\epsilon^{d-1}}\mu_1 \times \mu_2 \lbrace (x,y) \in E_1\times E_2 :
(x_d-y_d)t_1-\epsilon \leq x_1-y_1 \leq (x_d-y_d)t_1 + \epsilon, ...,
\\
&(x_d-y_d)t_{d-1}-\epsilon \leq x_{d-1}-y_{d-1} \leq (x_d-y_d)t_{d-1}
+ \epsilon \rbrace.
\end{align*}

Our plan is to show that $ \int \nu_\epsilon(t)$ is bounded above and below by a constant plus an error term. For, we write $\nu_\epsilon(t)$ as the sum of two terms. We shall prove that the integral of first term in $t_1, \dots, t_{d-1}$ is bounded above and below by two positive constants. Then we will show that the integral of the modulus of the second second term is bounded above by a sufficiently small positive constant plus a constant multiple of $\epsilon^{s-(d-1)}$. We will then conclude $1- \epsilon^{s-(d-1)} \lesssim \nu_{\epsilon}(t) \lesssim 1 + \epsilon^{s-(d-1)}$ for a subset of ${\left[ \frac{1}{2}, 1 \right]}^{d-1}$ of positive $(d-1)$-dimensional Lebesgue measure. As we note above, this implies that $\sigma({\mathcal D}(E))>0$. 

Let $\psi:\mathbb{R}\rightarrow \mathbb{R}$ be a smooth, even bump function, whose support is contained in the set $[-2,-\frac{1}{2}]\cup [\frac{1}{2},2]$ such that $\widehat{\psi}(0)=1$. We have 
\begin{align*}
\nu_\epsilon(t) &\approx \frac{1}{\epsilon^{d-1}}\int\int \psi\left(\frac{(x_1-y_1) -
t_1(x_d-y_d)}{\epsilon}\right)\psi\left(\frac{(x_2-y_2) -
t_2(x_d-y_d)}{\epsilon}\right)...\\
&\psi\left(\frac{(x_{d-1}-y_{d-1}) -
t_{d-1}(x_d-y_d)}{\epsilon}\right)d\mu_1(x)d\mu_2(y).
\end{align*}

By Fourier inversion, this quantity equals 

\begin{align*}
&\frac{1}{\epsilon^{d-1}} \int\int\int...\int
\widehat{\psi}(\lambda_1)\widehat{\psi}(\lambda_2)...\widehat{\psi}(\lambda_{d-1})
e^{\frac{2\pi i}{\epsilon}\lambda_1\left((x_1-y_1) -
t_1(x_d-y_d)\right)}e^{\frac{2\pi i}{\epsilon}\lambda_2\left((x_2-y_2) -
t_2(x_d-y_d)\right)} ...\\
&e^{\frac{2\pi i}{\epsilon}\lambda_{d-1}\left((x_{d-1}-y_{d-1}) -
t_{d-1}(x_d-y_d)\right)} d\lambda_1 d\lambda_2 ...
d\lambda_{d-1}d\mu_1(x)d\mu_2(y)\\
&=\frac{1}{\epsilon^{d-1}} \int\int...\int
\widehat{\psi}(\lambda_1)\widehat{\psi}(\lambda_2)...\widehat{\psi}(\lambda_{d-1})\widehat{\mu_1}\left(-\frac{\lambda_1}{\epsilon},
-\frac{\lambda_2}{\epsilon}, ...,\left(\frac{t_1\lambda_1 + t_2\lambda_2+...t_{d-1}\lambda_{d-1}}{\epsilon}\right)\right)
\\
&\overline{\widehat{\mu_2}}\left(-\frac{\lambda_1}{\epsilon},
-\frac{\lambda_2}{\epsilon}, ...,\left(\frac{t_1\lambda_1 + t_2\lambda_2+...t_{d-1}\lambda_{d-1}}{\epsilon}\right)\right)d\lambda_1
d\lambda_2...d\lambda_{d-1},\\
\end{align*}

and by changing variables $\lambda_j' = \frac{\lambda_j}{\epsilon}$, the latter is equal to

\begin{align*}
&\int\int...\int
\widehat{\psi}(\epsilon\lambda_1')\widehat{\psi}(\epsilon\lambda_2')...\widehat{\psi}(\epsilon\lambda_{d-1}')\widehat{\mu_1}\left(-\lambda_1',
-\lambda_2', ...,\left(t_1\lambda_1' + t_2\lambda_2'+...t_{d-1}\lambda_{d-1}'\right)\right)\\
&\overline{\widehat{\mu_2}}\left(-\lambda_1',
-\lambda_2', ...,\left(t_1\lambda_1' + t_2\lambda_2'+...t_{d-1}\lambda_{d-1}'\right)\right)d\lambda_1'
d\lambda_2'...d\lambda_{d-1}'.\\
\end{align*}

In what follows, we will refer to $\lambda_j'$ as  $\lambda_j$, to simplify the exposition. We also define

$$\lambda = \left( \lambda_1, \lambda_2, \dots, \lambda_{d-1} \right) \in \mathbb{R}^{d-1}.$$

For the sake of simplicity we define the function 

$$\Phi\left(\lambda\right)=\Phi\left(\lambda_1, \lambda_2, \dots,\lambda_{d-1}\right):=\psi(\lambda_1)\psi(\lambda_2)\dots\psi(\lambda_{d-1}).$$

An easy calculation shows that

$$\widehat{\Phi}(\lambda)=\widehat{\psi}(\lambda_1)\widehat{\psi}(\lambda_2)...\widehat{\psi}(\lambda_{d-1}),$$
which entails that $\widehat{\Phi}(0)=1$, since $\widehat{\psi}(0)=1$. Also, note that $\widehat{\Phi}$ is continuous. By the preceding discussion one can see that

\begin{align*}
\nu_\epsilon(t) &=\int\int\dots\int (\widehat{\Phi}(\epsilon \lambda)-\widehat{\Phi}(0)) \widehat{\mu_1}\left(-\lambda_1,
\lambda_2, ...,\left(t_1\lambda_1 + t_2\lambda_2+...t_{d-1}\lambda_{d-1}\right)\right)\\
&\overline{\widehat{\mu_2}}\left(\lambda_1,
\lambda_2, ...,\left(t_1\lambda_1 + t_2\lambda_2+...t_{d-1}\lambda_{d-1}\right)\right)d\lambda \\ 
&+\int\int\dots\int \widehat{\mu_1}\left(-\lambda_1,
\lambda_2, ...,\left(t_1\lambda_1 + t_2\lambda_2+...t_{d-1}\lambda_{d-1}\right)\right)\\
&\overline{\widehat{\mu_2}}\left(\lambda_1,
\lambda_2, ...,\left(t_1\lambda_1 + t_2\lambda_2+...t_{d-1}\lambda_{d-1}\right)\right)d\lambda \\ 
&=:P_\epsilon(t)+ M(t).
\end{align*}

We will prove that if $P_\epsilon(t)=Q_\epsilon(t)+R_\epsilon(t)$, where $Q_\epsilon(t)$ and $R_\epsilon(t)$ are error terms, for sufficiently small $\eta>0$,
$$\int M(t)dt \approx 1, \quad \int |Q_\epsilon(t)|dt \lesssim \eta I_s(\mu), \quad \text{and} \quad \int |P_\epsilon(t)|dt \lesssim \epsilon^{s-(d-1)} I_s(\mu).$$ 
Let $\delta$ be a small positive real number, to be chosen later. We handle the error term $P_\epsilon$ by splitting it into two integrals, where the domains of integration are $|\lambda | \leq \frac{\delta}{\epsilon}$ and $|\lambda |> \frac{\delta}{\epsilon}$ respectively.

\begin{equation*}
P_\epsilon = \int\int \dots\int_{|\lambda | \leq \frac{\delta}{\epsilon}} d\lambda_1
d\lambda_2...d\lambda_{d-1} + \int\int\dots\int_{|\lambda | > \frac{\delta}{\epsilon}} d\lambda_1
d\lambda_2...d\lambda_{d-1}\equiv Q_\epsilon(t)+R_\epsilon(t).
\end{equation*}

\vskip.125in

First we bound the $L^1$ norm of the quantity $R_\epsilon(t)$. Let $\psi_0 \in C^{\infty}_0(\mathbb{R})$ such that $supp(\psi_0) \subset [1/4,2]$ and $\psi_0=1$ in $[1/2,1]$.
\begin{align*}
& \int\limits_{[\frac{1}{2},1]^{d-1}} \psi_0(t_1) \dots \psi_0(t_{d-1}) |R_\epsilon(t)| dt_1...dt_{d-1}\\
& \leq \int\limits_{[\frac{1}{2},1]^{d-1}}\int\limits_{|\lambda | > \frac{\delta}{\epsilon}} |\widehat{\Phi}(\epsilon \lambda)-\widehat{\Phi}(0)| \left|\widehat{\mu_1}\left(-\lambda_1,
-\lambda_2, \dots,\left(t_1\lambda_1 + t_2\lambda_2+...t_{d-1}\lambda_{d-1}\right)\right)\right|\\
&\left|\overline{\widehat{\mu_2}}\left(-\lambda_1,
-\lambda_2, ...,\left(t_1\lambda_1 + t_2\lambda_2+...t_{d-1}\lambda_{d-1}\right)\right)\right| 
d\lambda_1 d\lambda_2...d\lambda_{d-1}dt_1...dt_{d-1}.
\end{align*}
By applying Cauchy-Schwarz, the square of the expression above is 
\begin{align*}
&\lesssim \int\limits_{[\frac{1}{2},1]^{d-1}}\int\limits_{|\lambda | > \frac{\delta}{\epsilon}}
|\widehat{\Phi}(\epsilon \lambda)-\widehat{\Phi}(0)| \psi_0(t_1) \psi_0(t_2)...\psi_0(t_{d-1}) \\
&\left | \widehat{\mu_1}\left(-\lambda_1,
-\lambda_2, ...,\left(t_1\lambda_1 + t_2\lambda_2+...t_{d-1}\lambda_{d-1}\right)\right)\right |^2 d\lambda_1 ...d\lambda_{d-1} dt_1...dt_{d-1}\\
&\int\limits_{[\frac{1}{2},1]^{d-1}}\int\limits_{|\lambda | > \frac{\delta}{\epsilon}} |\widehat{\Phi}(\epsilon \lambda)-\widehat{\Phi}(0)| \psi_0(t_1) \psi_0(t_2)...\psi_0(t_{d-1}) \\
&\left |\overline{\widehat{\mu_2}}\left(-\lambda_1,
-\lambda_2, ...,\left(t_1\lambda_1 + t_2\lambda_2+...t_{d-1}\lambda_{d-1}\right)\right)\right |^2 d\lambda_1...d\lambda_{d-1} dt_1...dt_{d-1}\\
& = A\cdot B,
\end{align*}
where $A$ is the first integral, and $B$ is the second. We will break each of these integrals up into $A_j$ and $B_j$, which are  integrals over subsets which make up the whole region of integration. Let us define
$$ \Lambda_j:= \{ \lambda \in \mathbb{R}^{d-1}: |\lambda| > \delta / \epsilon, \quad |\lambda_j| \approx |\lambda|\},$$

and notice that

$$A \leq \sum\limits_{j=1}^{d-1}\int\limits_{[\frac{1}{2},1]^{d-1}}\int\limits_{\Lambda_j} = \sum\limits_{j=1}^{d-1}A_j.
$$

We will now estimate the $A_j$, while the estimates of the corresponding $B_j$ are identical. Without loss of generality we may assume that $j=1$ and $\lambda_1>0$. The same proof works in the case $\lambda_1<0$ with minor modifications. For, we introduce the change of variables

$$
\tau_i=t_i,\quad \ell_j=-\lambda_j, \quad \ell_d=\sum_{k=1}^{d-1}t_k\lambda_k,
$$
where $2\leq i \leq d-1$ and $1 \leq j \leq d-1$. Notice that the Jacobian can be considered to have positive sign since otherwise we make a slightly different change of variables. Additionally, for the rest of the proof we keep denoting our new variables $\tau_i, \ell_j$ by $t_i, \lambda_j$ respectively, to stay consistent with our earlier notation. We, now, illustrate the estimation of the integral $A_1$, while the other pieces can be estimated in a similar fashion. 
\vskip.125in

By the definition of $\lambda_d$,
$$t_1 = \frac{\lambda_d + t_2\lambda_2+\dots t_{d-1}\lambda_{d-1}}{\lambda_1},$$
and we have that
\begin{align*}
A_1 &\leq \int_{[\frac{1}{2},1]^{d-2}}\psi_0(t_2)\dots\psi_0(t_{d-1})\int_{\Lambda_1} \int_{|\lambda_d|\lesssim |\lambda|}\psi_0\left( \frac{\lambda_d + t_2\lambda_2+\dots+t_{d-1}\lambda_{d-1}}{\lambda_1}\right)\\
&|\widehat{\Phi}(\epsilon \lambda)-\widehat{\Phi}(0)| \left|\widehat{\mu_1}(\lambda_1, \dots, \lambda_d)\right|^2 d\lambda_d \frac{d\lambda_1}{\lambda_1} d\lambda_2 \dots d\lambda_{d-1} dt_2 \dots dt_{d-1}.
\end{align*}

Define $\lambda' := (\lambda_1, \lambda_2, \dots, \lambda_{d-1}, \lambda_d)= (\lambda, \lambda_d) \in \mathbb{R}^d$ 
and one can easily deduce that 

$$\Lambda_1 \times \{\lambda_d: \lambda_d \lesssim |\lambda|\} \subset \{\lambda: |\lambda'| > \delta /\epsilon\}.$$

Hence, 

\begin{align*}
A_1 &\lesssim \int\int\dots\int_{|\lambda|>\frac{\delta}{\epsilon}} \left|\widehat{\mu_1}(\lambda_1, \dots, \lambda_d)\right|^2\frac{1}{|\lambda'|} d\lambda_1 \dots d\lambda_d \equiv \mathcal{I}_1.
\end{align*}
Now, we will use the energy integral bounds in polar coordinates. Recall that, since the Hausdorff dimension of $E$ is $s$, the energy integral of $\mu$ which is given by
$$
I_s({\mu}) :=\int_0^\infty\int_{S^{n-1}} \frac{\left|\widehat{\mu}(\xi)\right|^2}{|\xi|^{d-s}}d\xi=\int_0^\infty\int_{S^{n-1}} {\left|\widehat{\mu}(r\theta)\right|^2}\frac{r^{d-1}}{r^{d-s}}drd\theta,
$$

is finite. By rewriting the integral $\mathcal{I}_1$ in polar coordinates, we have

\begin{align*}
\mathcal{I}_1&= \int_{\delta / \epsilon}^\infty\int_{S^{n-1}}r^{d-s}\frac{\left|\widehat{\mu_1}(r\theta)\right|^2}{r}\frac{r^{d-1}}{r^{d-s}}drd\theta\\
&=\int_{\delta / \epsilon}^\infty\int_{S^{n-1}}\frac{\left|\widehat{\mu_1}(r\theta) \right|^2}{r^{s-(d-1)}}\frac{r^{d-1}}{r^{d-s}}drd\theta\\
&\leq \left(\frac{\epsilon}{\delta}\right)^{s-(d-1)} I_s(\mu),
\end{align*}

which concludes the proof of
\begin{equation}\label{errorterm}
\int \left|R_\epsilon(t) \right|dt \lesssim \epsilon^{s-(d-1)}.
\end{equation}

By the continuity of $\widehat{\Phi}$ for a sufficiently small positive constant $\eta$ we can find $\delta>0$ such that 
\begin{equation} \int |Q_\epsilon(t)| dt \lesssim \eta I_s(\mu).\label{mainerror}\end{equation}

It only remains to prove that for the main term,

\begin{equation}\label{mainterm}
\int M(t) dt \approx 1.
\end{equation}

The argument for the upper bound is similar to the one used in the proof of (\ref{mainerror}) but simpler since we do not use the continuity of $\widehat{\Phi}$. It only remains to prove the lower bound. To this end, by the same change of variables we made in the proof of (\ref{errorterm}) one can reduce case to proving a lower bound for the modulus of the integral

$$\int \int_{\Lambda_1'} \psi_0\left( \frac{\lambda_d + t_2\lambda_2+\dots+t_{d-1}\lambda_{d-1}}{\lambda_1}\right) \widehat{\mu_1}(\lambda_1, \dots, \lambda_d) \overline{\widehat{\mu_1}}(\lambda_1, \dots, \lambda_d)\frac{d\lambda_1 \dots d\lambda_d}{|\lambda'|}, $$

where $\Lambda_1':= \left\{ \lambda \in \mathbb{R}^{d-1}:|\lambda_1| \approx |\lambda|\right\}$. Moreover, since

$$\frac{\lambda_d + t_2\lambda_2+\dots+t_{d-1}\lambda_{d-1}}{\lambda_1} \in supp(\psi_0) \subset [\frac{1}{4},2]$$

we can reduce the case even more, by observing that the domain of integration can be considered as an appropriate $d$-dimensional sector $\mathcal{S}$. By readjusting the constants we may assume that the sector is of angle $\pi/3$. Then, if we let $\Theta$ to be the rotation by $\pi/3$ we can write

$$\sum_{i=0}^5 \int_{\Theta^i \mathcal{S}} \widehat{\mu_1}(\lambda_1, \dots, \lambda_d) \overline{\widehat{\mu_2}}(\lambda_1, \dots, \lambda_d)\frac{d\lambda_1 \dots d\lambda_d}{|\lambda'|}= \int\widehat{\mu_1}(\lambda_1, \dots, \lambda_d) \overline{\widehat{\mu_2}}(\lambda_1, \dots, \lambda_d)\frac{d\lambda_1 \dots d\lambda_d}{|\lambda'|}.$$

Since $E_1, E_2 \subset [0,1]^d$ are separated by the stopping time argument employed previously, one can assume that

$$\int \widehat{\mu_1}(\lambda_1, \dots, \lambda_d) \overline{\widehat{\mu_2}}(\lambda_1, \dots, \lambda_d)\frac{d\lambda_1 \dots d\lambda_d}{|\lambda'|} \approx \int \int \frac{d\mu_1(x) d\mu_2(y)}{|x-y|} \approx 1.$$

Thus, there exists $i_0 \in \{0,1,\dots,5\}$ such that 

$$\int_{\Theta^{i_0} \mathcal{S}} \widehat{\mu_1}(\lambda_1, \dots, \lambda_d) \overline{\widehat{\mu_2}}(\lambda_1, \dots, \lambda_d)\frac{d\lambda_1 \dots d\lambda_d}{|\lambda'|} \gtrsim 1,$$

and if $\Theta^{i_0}:= \Theta_0$, by a change of variables we obtain

\begin{equation}
\int_\mathcal{S} \widehat{\mu_1}(\Theta_0^{-1} \lambda') \overline{\widehat{\mu_2}}(\Theta_0^{-1} \lambda')\frac{d\lambda'}{|\lambda'|} \gtrsim 1.
\label{energytheta}\end{equation}

Let us define, now, the measure $\mu_j^{\Theta_0}$ on the set $E_j^{\Theta_0}:= \{\Theta_0^{-1}x: x\in E_j\},$ for $j=1,2,$ as follows:

$$\int f(x) d\mu^{\Theta_0}_j(x)=\int f(\Theta^{-1}_0 x) d\mu_j(x).$$

The latter, for $f(x)=e^{-2\pi i x \cdot y}$, shows that 

$$\widehat{\mu^{\Theta_0}_j}(y)=\widehat{\mu_j}(\Theta_0^{-1}y),$$

and thus, (\ref{energytheta}) may be written as

\begin{equation}
\int_\mathcal{S} \widehat{\mu^{\Theta_0}_1}( \lambda') \overline{\widehat{\mu^{\Theta_0}_2}( \lambda')}\frac{d\lambda'}{|\lambda'|} \gtrsim 1.
\label{energytheta2}\end{equation}

The estimate (\ref{energytheta2}) provides us with the lower bound we were looking for, for the measures $\mu^{\Theta_0}_j$ instead of $\mu_j$. This allows us to prove that the set of directions determined by the sets $E_1^{\Theta_0}$ and $E_2^{\Theta_0}$ has positive Lebesgue measure.To this end, we define the measure $\nu^{\Theta_0}$ with support, without loss of generality, the closed cube $[1/2,1]^{d-1}$, as follows:

$$ \int g(t) d\nu^{\Theta_0}(t) = \int g\left( \frac{x_1-y_1}{x_d-y_d}, \frac{x_2-y_2}{x_d-y_d},\dots, \frac{x_{d-1}-y_{d-1}}{x_d-y_d}\right) d\mu^{\Theta_0}_1(x)d\mu^{\Theta_0}_2(y)$$

and the quantity $\nu_\epsilon^{\Theta_0}(t)$ to be equal to

$$\frac{1}{\epsilon^{d-1}}\mu^{\Theta_0} \times \mu^{\Theta_0} \left\lbrace (x,y) \in E^{\Theta_0}_1\times E^{\Theta_0}_2 : t_1-\epsilon \leq \frac{x_1-y_1}{x_d-y_d} \leq t_1 + \epsilon, ..., t_{d-1}-\epsilon \leq \frac{x_{d-1}-y_{d-1}}{x_d-y_d} \leq t_{d-1} + \epsilon \right\rbrace,$$

or else, by a stopping time argument, comparable to

$$\frac{1}{\epsilon^{d-1}} \mu^{\Theta_0}_1 \times \mu^{\Theta_0}_2 \left\lbrace (x,y) \in E^{\Theta_0}_1\times E^{\Theta_0}_2 :
|x_1-y_1-(x_d-y_d)t_1| \leq \epsilon,\dots,|x_{d-1}-y_{d-1}-(x_d-y_d)t_{d-1}| \leq \epsilon\right\rbrace.$$

We follow $mutatis$ $mutandi$ the preceding argument for $\nu_\epsilon(t)$, and in view of (\ref{energytheta2}) we prove that

$$ \int M^{\Theta_0}(t) dt \approx 1, \quad  \int |Q_\epsilon^{\Theta_0}(t)| dt \lesssim \eta \quad \text{and} \quad \int |R_\epsilon^{\Theta_0}(t)| dt \lesssim \epsilon^{s-(d-1)}.$$

Therefore, since $\eta>0$ is sufficiently small

\begin{equation}1 -\epsilon^{s-(d-1)} \lesssim \int \nu_\epsilon^{\Theta_0}(t)dt \lesssim 1+ \epsilon^{s-(d-1)}.\label{nubounds}\end{equation}

By Lebesgue's decomposition 

$$d\nu^{\Theta_0} = d\nu^{\Theta_0}_{ac} + d\nu^{\Theta_0}_s,$$

where $\nu^{\Theta_0}_{ac} \ll \mathcal{L}^{d-1}$ and $\nu^{\Theta_0}_s \perp \mathcal{L}^{d-1}$. Moreover, if $f$ is the Radon-Nikodym derivative of $\nu^{\Theta_0}_{ac}$ with respect to $\mathcal{L}^{d-1}$ we can write $d\nu^{\Theta_0}_{ac}=f d\mathcal{L}^{d-1}$. Since $\nu^{\Theta_0}$ is a finite Borel measure, 

$$ \lim_{\epsilon \to 0}\frac{\nu^{\Theta_0}(B(t,\epsilon))}{\mathcal{L}^{d-1}(B(t,\epsilon))}=f(t),$$

for  $\mathcal{L}^{d-1}$-a.e. $t$. By (\ref{nubounds}) and Fatou's lemma

\begin{equation}\int \liminf_{\epsilon \to 0}  \frac{\nu^{\Theta_0}(B(t,\epsilon))}{\mathcal{L}^{d-1}(B(t,\epsilon))} dt \leq \int \liminf_{\epsilon \to 0}\frac{\nu^{\Theta_0}(B(t,\epsilon))}{\mathcal{L}^{d-1}(B(t,\epsilon))} dt = \int f(t) dt \lesssim 1. \label{radnikupper}\end{equation}

where $B(t,\epsilon)$ is the $(d-1)$-dimensional ball with center at $t$ and of radius $\epsilon$. Note that 
$$f <\infty, \quad \mathcal{L}^{d-1}\text{-a.e.}$$

Let us now define for an appropriate constant $C_0>0$,

$$ \nu^{\Theta_0}_*(t) :=  \sup_{\epsilon< C_0} \frac{\nu^{\Theta_0}(B(t,\epsilon))}{\mathcal{L}^{d-1}(B(t,\epsilon))}.$$

Observe that $\nu^{\Theta_0}_*(t)$ provides us with a dominating function for $\nu^{\Theta_0}_{\epsilon}(t)$ and by our previous arguments for the upper bounds 

$$\int \nu^{\Theta_0}_*(t) dt \lesssim 1.$$

Fatou's lemma for $\nu^{\Theta_0}_*(t)-\nu^{\Theta_0}_{\epsilon}(t)$ and (\ref{nubounds}) show that

$$1 \lesssim \limsup_{\epsilon \to 0} \int \frac{\nu^{\Theta_0}(B(t,\epsilon))}{\mathcal{L}^{d-1}(B(t,\epsilon))}dt \leq \int \limsup_{\epsilon \to 0} \frac{\nu^{\Theta_0}(B(t,\epsilon))}{\mathcal{L}^{d-1}(B(t,\epsilon))}dt = \int f(t) dt,$$

which in view of (\ref{radnikupper}) gives

\begin{equation}\int f(t) dt \approx 1. \label{radnikapprox}\end{equation}

Let $\mathcal{F}$ be the support of $\nu^{\Theta_0}_{ac}$ and note that 

$$\mathcal{L}^{d-1}(\mathcal{F})=\mathcal{L}^{d-1}([1/2,1]^{d-1}), \quad \nu^{\Theta_0}\lfloor \mathcal{F}=\nu^{\Theta_0}_{ac} \quad \text{and} \quad \nu^{\Theta_0}\lfloor\mathcal{F} \ll \mathcal{L}^{d-1}.$$

 Thus for every Borel set $B$ in $\mathcal{F}$

$$ \nu^{\Theta_0}\lfloor \mathcal{F} (B) = \int_B f(t) dt$$

which setting $B=\mathcal{F}$ in conjunction with (\ref{radnikapprox}) implies that the set of directions determined by $E^{\Theta_0}_1$ and $E^{\Theta_0}_2$ contains a set of positive Lebesgue measure and since, by definition, $E^{\Theta_0}_1$ and $E^{\Theta_0}_2$ determine the same number of directions with $E_1$ and $E_2$ the proof of Theorem \ref{main} is concluded.

\section{Some connections between continuous and discrete aspects of the problem at hand} 
\label{discrete} 

\vskip.125in 

In this section we appeal to a conversion mechanism developed in \cite{IL05}, \cite{HI05}, and \cite{IRU10}, to deduce a Pach-Pinchasi-Sharir type result from Theorem \ref{main}. In the aforementioned papers, the conversion mechanism was used in the context of distance sets. However, as we shall see below, the idea is quite flexible and lends itself to a variety of applications. 

\begin{definition} \label{adaptablemama} Let $P$ be a set of $n$ points contained in ${[0,1]}^d$, $d \ge 2$. Define the measure
$$ d \mu^s_P(x)=n^{-1} \cdot n^{\frac{d}{s}} \cdot \sum_{p \in P} \chi_{B_{n^{-\frac{1}{s}}}(p)}(x)dx,$$ where $\chi_{B_{n^{-\frac{1}{s}}}(p)}(x)$ is the characteristic function of the ball of radius $n^{-\frac{1}{s}}$ centered at $p$. 

We say that $P$ is $s$-adaptable if $P$ is $n^{-\frac{1}{s}}$-separated and 

$$ I_s(\mu_P)=\int \int {|x-y|}^{-s} d\mu^s_P(x) d\mu^s_P(y)<\infty.$$ 

\vskip.125in 

This is equivalent to the statement 

$$ n^{-2} \sum_{p \not=p' \in P} {|p-p'|}^{-s} \lesssim 1.$$ 

\end{definition}

\vskip.125in 

To put it simply, $s$-adaptability means that a discrete point set $P$ can be thickened into a set which is uniformly $s$-dimensional in the sense that its energy integral of order $s$ is finite. Unfortunately, it is shown in \cite{IRU10} that there exist finite point sets which are not $s$-adaptable for certain ranges of the parameter $s$. The point is that the notion of Hausdorff dimension is much more subtle than the simple ``size" estimate. This turns out to be a serious obstruction to efforts to convert ``continuous" results into ``discrete analogs". 

\vskip.125in 

The first main result of this section is the following. 

\begin{theorem} \label{continuoustodiscrete} Suppose that for arbitrarily small $\epsilon>0$, $P$ is 
a $(d-1+\epsilon)$-adaptable set in ${[0,1]}^d$, $d \ge 2$, consisting of $n$ points. Then 

\begin{equation} \label{discreteorgasm} \# {\mathcal D}(P) \gtrapprox n. \end{equation} 

Moreover, there exists a subset, ${\mathcal D}'(P) \subset {\mathcal D}(P)$, such that $\# {\mathcal D}'(P) \ge \frac{1}{2} \# {\mathcal D}(P)$ and the elements in ${\mathcal D}'(P)$ are $n^{-\frac{d-1}{d-1+\epsilon}}$-separated.
\end{theorem} 

Observe that in dimensions two and three, this result is much weaker than what is known, as we note in the introduction above. Another weakness of this result is that it only holds for $s$-adaptable sets. However, in dimensions four and higher, Theorem \ref{continuoustodiscrete} appears to give a new result in the discrete setting. 

To prove Theorem \ref{continuoustodiscrete} thicken each point of $P$ by $n^{-\frac{1}{s}}$, where $s>d-1$. Let $E_P$ denote the resulting set. Then 

$$ \sigma({\mathcal D}(E_P)) \lesssim n^{-\frac{d-1}{s}} \cdot \# {\mathcal D}(P).$$ 

\vskip.125in

By the adaptability assumption and the proof of Theorem \ref{main}, we see that 

$$ \# {\mathcal D}(P) \gtrapprox n^{\frac{d-1}{s}},$$ establishing (\ref{discreteorgasm}) in view of the fact that we may take $s$ arbitrarily close to $d-1$. Note that since Theorem \ref{main} does not hold for $s=d-1$, we cannot replace (\ref{discreteorgasm}) by 

\begin{equation} \label{multiplediscreteorgasms} \# {\mathcal D}(P) \ge C \# P. \end{equation} 

To see that there exists a subset of the direction set which is $n^{-\frac{d-1}{d-1+\epsilon}}$-separated, we recall that $\sigma({\mathcal D}(E_P))>0$, so we can break ${\mathcal D}(E_P)$ up into pieces with Lebesgue measure $n^{-\frac{d-1}{d-1+\epsilon}}$, each of which contains a representative from ${\mathcal D}(P)$. Then by a simple pigeon-hole argument, we see that at least $\displaystyle \frac{1}{2^{d-1}}\#{\mathcal D}(P)$ of these must be separated.

\begin{remark}
One should take note that the separation statement is relatively unique to continuous techniques. Most of the standard discrete results say nothing about the separation of the distinct elements in a given set. For example, in \cite{PPS04}, there is a sharp lower bound on the number of distinct directions determined by a set of points in $\mathbb{R}^3$, but there are no guarantees on the separation or distribution of these directions on $S^2$.
\end{remark}


\newpage

\begin{thebibliography}{8}

\bibitem{BMP00} P. Brass, W. Moser, and J Pach, {Research Problems in Discrete Geometry}, Springer (2000). 

\bibitem{B86} J. Bourgain, {\it A Szemerdi type theorem for sets of positive density}, Israel J. Math. \textbf{54} (1986), no. 3, 307-331.

\bibitem{B94} J. Bourgain, {\it Hausdorff dimension and distance sets} Israel. J. Math. \textbf{87} (1994), 193-201. 

\bibitem{CEHIT10} D. Covert, B. Erdogan, D. Hart, A. Iosevich, and K. Taylor {\it Finite point configurations, uniform distribution, intersections of fractals, and number theoretic consequences}, (in preparation), (2010). 

\bibitem{Erd05} B. Erdo\u{g}an {\it A bilinear Fourier extension theorem and applications to the distance set problem} IMRN (2006).

\bibitem{Fal86} K. J. Falconer {\it On the Hausdorff dimensions of distance sets} Mathematika \textbf{32} (1986) 206-212.

\bibitem{Falc86} K. J. Falconer, {\it The geometry of fractal sets}, Cambridge Tracts in Mathematics, \textbf{85} Cambridge University Press, Cambridge, (1986). 

\bibitem{FKW90} H. Furstenberg, Y. Katznelson, and B. Weiss, {\it Ergodic theory and configurations in sets of positive density} Mathematics of Ramsey theory, 184-198, Algorithms Combin., 5, Springer, Berlin, (1990).

\bibitem{GV03} J. Garnett and J. Verdera {\it Analytic capacity, bilipschitz maps and Cantor sets}, Math. Res. Lett. \textbf{10} (2003), no. 4, 515-522. 

\bibitem{HI05} S. Hofmann and A. Iosevich, A {\it Circular averages and Falconer/Erd\H os distance conjecture in the plane for random metrics}, Proc. Amer. Math. Soc. \textbf{133} (2005), no. 1, 133-143.

\bibitem{HIKR10} D. Hart, A. Iosevich, D. Koh, and M. Rudnev, {\it Averages over hyperplanes, sum-product theory in vector spaces over finite fields and the Erd\H{o}s-Falconer distance conjecture}, Transaction of the AMS, (in press), (2010). 

\bibitem{IL05} A. Iosevich, I. Laba, {\it$K$-distance sets, Falconer conjecture, and discrete analogs}, Integers \textbf{5} (2005), no. 2. 

\bibitem{IM10} A. Iosevich, H. Morgan and J. Pakianathan, {\it Sets of directions determined by subsets of vector spaces over finite fields}, (preprint), (2010). 

\bibitem{IRU10} A. Iosevich, M. Rudnev, and I. Uriarte-Tuero, {\it Theory of dimension for large discrete sets and applications}, (arXiv:0707.1322). 

\bibitem{IS10} A. Iosevich and S. Senger, {\it Sharpness of Falconer's estimate in continuous and arithmetic settings, geometric incidence theorems and distribution of lattice points in convex domains}, (submitted for publication) (2010). 

\bibitem{Landau27} E. Landau, {\it Zahlentheorie}, Chelsie Publishing Company, (1927). 

\bibitem{Ma02} J. Matousek, {\it Lectures on Discrete Geometry,} Graduate Texts in Mathematics, Springer \textbf{202} (2002).

\bibitem{Mat87} P. Mattila {\it Spherical averages of Fourier transforms of measures with finite energy: dimensions of intersections and distance sets} Mathematika, \textbf{34} (1987),  207-228.

\bibitem{M95} P. Mattila, {\it Geometry of sets and measures in Euclidean spaces}, Cambridge University Press, \text{volume} 44, (1995). 

\bibitem{OS10} T. Orponen and T. Sahlsten, {\it Radial projections of rectifiable sets}, (http://arxiv.org/pdf/1101.2587), Annales Academi{\ae} Scientiarum Fennic{\ae} Mathematic (accepted for publication) (2010). 

\bibitem{P05} J. Pach {\it Directions in combinatorial geometry}, Jahresber. Deutsch. Math.-Verein. 
\textbf{107} (2005), no. 4, 215-225. 

\bibitem{PS04} J. Pach and M. Sharir {\it Geometric incidences}, Towards a theory of geometric graphs, 185-223, Contemp. Math., \textbf{342}, Amer. Math. Soc., Providence, RI, (2004). 

\bibitem{PPS04} J. Pach, R. Pinchasi, M. Sharir {\it  On the number of directions determined by a three-dimensional points set}, J. Combin. Theory Ser. A \textbf{108} (2004), no. 1, 1-16.

\bibitem{PPS07} J. Pach, R. Pinchasi, M. Sharir {\it Solution of Scott's problem on the number of directions determined by a point set in 3-space}, Discrete Comput. Geom. \textbf{38} (2007), no. 2, 399-441.

\bibitem{SiSo06} K. Simon and B. Solomyak {\it Visibility for self-similar sets of dimension one in the plane} Real Anal. Exchange \textbf{32} (2006/07), no. 1, 67-78.

\bibitem{Sz97} L. Sz\'{e}kely, {\it A. Crossing numbers and hard Erd\H os problems in discrete geometry} Combin. Probab. Comput. \textbf{6} (1997), 353-358.

\bibitem{W99} T. Wolff, {\it Decay of circular means of Fourier transforms of measures}, International Mathematics Research Notices \textbf{10} (1999) 547-567.

\bibitem{W99II} T. Wolff, {Recent work connected with the Kakeya problem}, Prospects in Mathematics, H. Rossi, ed., American Mathematical Society, (1999).

\bibitem{Z06} T. Ziegler, {\it Nilfactors of ${\Bbb R}^d$ actions and configurations in sets of positive upper density in ${\Bbb R}^m$}, J. Anal. Math. \textbf{99}, pp. 249-266 (2006). 



\end{thebibliography}
\end{document}